\documentclass[11pt]{amsart}%
\usepackage{amsmath}%
\usepackage{amssymb}%
\usepackage{amscd}%
\usepackage{color}%
\usepackage{mathrsfs}%
\usepackage[all]{xy}%
\def\ol#1{\overline{#1}}
\def\wh#1{\widehat{#1}}
\def\wt#1{\widetilde{#1}}
\def\ul#1{\underline{#1}}
\theoremstyle{plain}
    \newtheorem{theorem}{Theorem}[section]
    \newtheorem{conjecture}[theorem]{Conjecture}
    \newtheorem{proposition}[theorem]{Proposition}
    \newtheorem{lemma}[theorem]{Lemma}
    \newtheorem{corollary}[theorem]{Corollary}
\theoremstyle{definition}
    \newtheorem{definition}[theorem]{Definition}

    \newtheorem{remark}[theorem]{Remark}
\def\Alphabet{A,B,C,D,E,F,G,H,I,J,K,L,M,N,O,P,Q,R,S,T,U,V,W,X,Y,Z}
\def\alphabet{a,b,c,d,e,f,g,h,i,j,k,l,m,n,o,p,q,r,s,t,u,v,w,x,y,z}
\def\endpiece{xxx}
\def\makeAlphabet[#1]{\expandafter\makeA#1,xxx,}
\def\makealphabet[#1]{\expandafter\makea#1,xxx,}
\def\makeA#1,{\def\temp{#1}\ifx\temp\endpiece\else%
\mkbb{#1}\mkfrak{#1}\mkbf{#1}\mkcal{#1}\mkscr{#1}\expandafter\makeA\fi}%
\def\makea#1,{\def\temp{#1}\ifx\temp\endpiece\else\mkfrak{#1}\mkbf{#1}\expandafter\makea\fi}%
\def\mkbb#1{\expandafter\def\csname bb#1\endcsname{\mathbb{#1}}}
\def\mkfrak#1{\expandafter\def\csname fr#1\endcsname{\mathfrak{#1}}}
\def\mkbf#1{\expandafter\def\csname b#1\endcsname{\mathbf{#1}}}
\def\mkcal#1{\expandafter\def\csname c#1\endcsname{\mathcal{#1}}}
\def\mkscr#1{\expandafter\def\csname s#1\endcsname{\mathscr{#1}}}
\def\makeop[#1]{\xmakeop#1,xxx,}
\def\mkop#1{\expandafter\def\csname #1\endcsname{{\mathrm{#1}}}} %
\def\xmakeop#1,{\def\temp{#1}\ifx\temp\endpiece\else\mkop{#1}\expandafter\xmakeop\fi}%
\makeAlphabet[\Alphabet]
\makealphabet[\alphabet]
\makeop[mot,dR,syn,Gal,ur]
\makeop[Hom,Ext,Sym,Coker]
\makeop[Eis,pol,Tr,Frob,ord,Fil]
\makeop[Spec,Spf,Spm,cyc,tg,cris,max,rig,pr]

\def\thepi{\boldsymbol{\pi}}
\def\levelstructure{\nu}
\def\cohomologyclass{\alpha}
\def\thealpha{u}
\def\thebeta{v}
\def\pair#1{\langle #1 \rangle}
\def\dmu{\mu}
\def\dmud{\mu'}
\def\dwtmu{\wt\mu}

\begin{document}
\title[$p$-adic Beilinson conjecture]{$p$-adic Beilinson conjecture for ordinary Hecke motives associated to imaginary quadratic fields}
\author{Kenichi Bannai}
\address{Department of Mathematics, Keio University, 3-14-1 Hiyoshi, Kouhoku-ku, Yokohama, JAPAN}
\email{bannai@math.keio.ac.jp}
\author{Guido Kings}
\address{NWF-I Mathematik, Universit\"at Regensburg, 93040 Regensburg, Germany}
\email{guido.kings@mathematik.uni-regensburg.de}
\date{\today }
\thanks{This work was supported by KAKENHI 21674001.}

\maketitle

%
%
%
\section{Introduction}
%
%
%

The purpose of this article is to give an overview of  the series of papers
\cite{BK1} \cite{BK2} concerning the $p$-adic Beilinson conjecture of motives associated to Hecke characters of 
an imaginary quadratic field $K$, for a prime $p$ which splits in $K$. The $p$-adic $L$-function
for such $p$ interpolating \textit{critical values} of $L$-functions of Hecke characters associated to imaginary 
quadratic fields was first constructed by Vishik and Manin \cite{VM}, and a different  construction using $p$-adic Eisenstein 
series was given by Katz \cite{Ka1}.  
The $p$-adic Beilinson conjecture, as formulated by Perrin-Riou in \cite{PR},
gives a precise conjecture concerning the \textit{non-critical values} of $p$-adic $L$-functions associated to general motives.
The purpose of our research is to investigate the interpolation property at non-critical points of the $p$-adic $L$-function 
constructed by Vishik-Manin and Katz.  

For simplicity, we assume in this article that the imaginary quadratic field $K$ has
class number \textit{one} and that the Hecke character $\psi$ we consider corresponds to an elliptic curve with complex multiplication
defined over $\bbQ$.  Let $a$ be an integer $>0$.
The main theorem of this article (Theorem \ref{thm: main}) is a proof of the $p$-adic Beilinson conjecture
for $\psi^a$ (see Conjecture \ref{conj: PR}),  when the prime $p \geq 5$ is an \textit{ordinary} prime.
The authors would like to thank the organizers Takashi Ichikawa, Masanari Kida and  Takao Yamazaki for
the opportunity to present our research at the RIMS ``Algebraic Number Theory and Related Topics 2009" 
conference.

%
%
%
\section{The $p$-adic Beilinson conjecture}\label{section: p Beilinson}
%
%
%

Assume that $K$ is an imaginary quadratic field of \textit{class number one}.
Let $E$ be an elliptic curve defined over $\bbQ$. 
We assume in addition that $E$ has complex multiplication by the ring of integers $\cO_K$ of $K$.
We let $\psi := \psi_{E/K}$ be the Grossencharacter of $K$ associated to $E_K := E \otimes_\bbQ K$ by the theory of complex 
multiplication, and we denote by $\frf$ the conductor of $\psi$.

We let $M(\psi)$ be the motive over $K$ with coefficients in $K$ associated to the Grossencharacter $\psi$.
Then we have $M(\psi) = H^1(E_K)$, where $H^1(E_K)$ is the motive associated to $E_K$. 
The Hasse-Weil $L$-function of $M(\psi)$ is a function with values in $K \otimes_\bbQ \bbC$ given by
$$
	L(M(\psi), s) = ( L(\psi_\tau,s))_{\tau : K \hookrightarrow \bbC},
$$
where $\tau: K \hookrightarrow \bbC$ are the embeddings of the coefficient $K$ of $M(\psi)$ into $\bbC$
and $L(\psi_\tau,s)$ is the Hecke $L$-function 
$$
	L(\psi_\tau,s) = \prod_{(\frq,\frf)=1} \left(1 - \frac{\psi_\tau(\frq)}{N\frq^s} \right)^{-1}
$$
associated to the character
$
	\psi_\tau: \bbA^\times_K \xrightarrow\psi K \overset{\tau}\hookrightarrow \bbC.
$
Here, the product is over the prime ideals $\frq$ of $K$ which are prime to $\frf$.

For integers $a>0$ and $n$, we let $M^a = M(\psi^{a}) := M(\psi)^{\otimes_K a}$, which is a motive over $K$ with coefficients in $K$.
Then the Hasse-Weil $L$-function $L(M^a, s)$ is given by the Hecke $L$-function
$$
	L(M^a, s) = (L(\psi_\tau^a, s))_{\tau: K \hookrightarrow \bbC}
$$
with values in $K \otimes_\bbQ\bbC$.
We let $M^a_B$ be the Betti realization of $M^a$, which is a $K$-vector space 
of dimension \textit{one}.  We fix a $K$-basis $\omega^a_B$ of $M^a_B$.
The de Rham realization $M^a_\dR(n)$ of $M^a(n)$  is the rank one $K \otimes_\bbQ K$-module
$$
	M^a_\dR(n) = K  \omega^{n-a,n} \bigoplus K \omega^{n,n-a},
$$
with Hodge filtration given by
$$
	F^m M^a_\dR(n) =
	\begin{cases}
		 M^a_\dR(n) &  m \leq -n \\
		 K\omega^{n,n-a}  &   -n < m \leq a-n \\
		 0  & \text{otherwise}. 
	\end{cases}
$$

In what follows, we consider the case when $n>a$, which implies in particular 
that our  motive is \textit{non-critical}.  
We have in this case $F^0 M^a_\dR(n) = 0$.  The tangent space of our motive is given by
$$
	t^a_{n} := M^a_\dR(n)/ F^0 M^a_\dR(n) \cong M^a_\dR(n),
$$
which is again a $K \otimes_\bbQ K$-module of rank one.   
Note that $\omega^a_{\tg,n} : = \omega^{n-a,n} + \omega^{n,n-a}$ gives a basis of 
$t^a_n$ as a $K\otimes_\bbQ K$-module.

We denote by $V^a_\infty(n)$ the $\bbR$-Hodge realization of  $M^a(n)$.
The Beilinson-Deligne cohomology 
$
	H^1_\sD(K \otimes_\bbQ \bbR,  V^a_\infty(n) )
$
is given as the cokernel of the natural inclusion
$$
	M^a_B(n) \otimes_\bbQ\bbR \rightarrow t^a_n \otimes_\bbQ \bbR.
$$
The Beilinson regulator map gives a homomorphism
\begin{equation}\label{eq: Beilinson regulator}
	r_\infty : H^1_\mot(K, M^a(n)) \rightarrow H^1_\sD( K \otimes_\bbQ \bbR, V^a_\infty(n)),
\end{equation}
from the motivic cohomology $H^1_\mot(K, M^a(n))$ of $K$ with coefficients in $M^a(n)$
to $H^1_\sD( K \otimes_\bbQ \bbR, V^a_\infty(n))$.  Then $r_\infty \otimes_\bbQ \bbR$
is known to be surjective and is conjectured to be an isomorphism.
We let $c^a_n$ be an element of $H^1_\mot(K, M^a(n))$ such that $r_\infty(c^a_n)$ generates
$H^1_\sD( K_\frp, V^a_\infty(n))$ as a $K_\infty := K \otimes_\bbQ \bbR$-module. 
We define the complex period $\Omega_\infty(n)$ of $M^a(n)$ to be the determinant of the
exact sequence
\begin{equation}\label{eq: complex determinant}
	0 \rightarrow M^a_B(n) \otimes_\bbQ \bbR \rightarrow t^a_n \otimes_\bbQ \bbR
	\rightarrow H^1_\sD(K \otimes_\bbQ \bbR,  V^a_\infty(n) ) \rightarrow 0
\end{equation}
for the basis $r_\infty(c^a_n)$, $\omega_{\tg,n}^a$, and $\omega^a_B$.  The complex period
is an element in $K_\infty$ and is independent of the choice of the basis up to multiplication 
by an element in $K^\times$.  The value $L(M^a,n)$ is in $K \otimes_\bbQ \bbR$, and
the weak Beilinson conjecture for $M^a(n)$ as proved by Deninger  \cite{Den} gives the following
(see Theorem \ref{thm: Deninger} and Corollary \ref{cor: Deninger} for the precise statement.)

\begin{theorem}\label{thm: Deninger one}
	The value
	$$
			\frac{L(M^a, n)}{\Omega_\infty(n)}
	$$
	is an element in $K^\times$.
\end{theorem}

For any prime $p$, the \'etale realization  $V^a_p(n)$ of our motive is a $K \otimes_\bbQ \bbQ_p$-vector space with continuous action of 
$\Gal(\ol K/K)$.  We fix a prime $p \geq 5$ relatively prime to $\frf$  such that $E$ has good \textit{ordinary} reduction at $p$. 
In this case, the ideal generated by $p$ splits as $(p) = \frp\frp^*$ in $K$.  
We fix a prime ideal $\frp$ of $K$ above $p$.  
Then the Bloch-Kato exponential map gives an isomorphism
\begin{equation}\label{eq: exponential map}
	\exp_p : t^a_n \otimes_K K_\frp \xrightarrow\cong H^1_f( K_\frp, V^a_p(n)),
\end{equation}
and the inverse of this isomorphism is denoted by $\log_p$.
The $p$-adic \'etale regulator map gives a homomorphism
\begin{equation}\label{eq: p-adic regulator}
	r_p : H^1_\mot(K, M^a(n)) \rightarrow H^1_f( K_\frp, V^a_p(n)),
\end{equation}
and the map $r_p \otimes \bbQ_p$ is conjectured to be an isomorphism.  Assuming that the $p$-adic regulator map $r_p$ is 
\textit{injective}, we define the $p$-adic period $\Omega_p(n)$ of $M^a(n)$ to be the determinant of the map 
$\log_p$ for the basis $r_p(c^a_n)$ and $\omega_{\tg,n}^a$.  
In other words, $\Omega_p(n)$ is an element in $K_p := K \otimes_\bbQ K_\frp \cong K_\frp \bigoplus K_{\frp^*}$ satisfying
\begin{equation}\label{eq: p-adic period}
	\log_p \circ\,\, r_{p}(c^a_n) = \Omega_p(n) \,\omega_{\tg,n}^a.
\end{equation}
The $p$-adic period $\Omega_p(n)$ is independent of the choice of basis up to multiplication by an element in $K^\times$.

\begin{remark}
	We need to assume the injectivity of the $p$-adic regulator $r_p$ to insure that the $p$-adic period $\Omega_p(n)$ is non-zero.
	Kato has proved in \cite{Kato} 15.15 the weak Leopoldt conjecture for any Hecke character of $K$.
	Hence by a result of Jannsen (\cite{Jan}, Lemma 8),  we may then conclude that 
	$$
		H^2(\cO_K[1/p\frf_a],V^a_p(n))=0
	$$
	for almost all $n$. This implies that $r_p$ is injective for such $n$.
\end{remark}

The $p$-adic Beilinson conjecture as formulated by Perrin-Riou (see \cite{Col} Conjecture 2.7) 
specialized to our setting is given as follows.

\begin{conjecture}\label{conj: PR}
	Let $a$ be an integer $>0$. Then there exists a $p$-adic pseudo-measure $\dmu^a$ on $\bbZ_p$ with values in $K_p$
	such that the value
	$$
		L_p(\psi^a  \otimes \chi^n_\cyc) := \int_{\bbZ^\times_p} w^n \dmu^a(w)
	$$
	in $K_p$ for any integer $n  > a$ satisfies
	$$
		\frac{L_p(\psi^a  \otimes \chi^n_\cyc)}{\Omega_p(n)}
		 = \left(1 - \frac{\psi(\frp)^a}{p^n} \right) \left(1 - \frac{\ol{\psi}(\frp^*)^a}{p^{a+1-n}} \right)  \frac{\Gamma(n) L(\psi^a, n)}{\Omega_\infty(n)},
	$$
	where $\frp$ is a fixed prime in $K$ above $p$. 
\end{conjecture}
If $\frf_a\neq(1)$ for the conductor $\frf_a$ of $\psi^a$, then $\dmu^a$ should in fact be a $p$-adic measure.
Note that the dependence of the pseudo-measure on the choices of the basis $\omega^a_{\tg,n}$ and $c^a_n$ cancel,
where as the  pseudo-measure depends on the choice of the basis $\omega^a_B$.
The main goal of our research is to prove that the $p$-adic measure constructed by Vishik-Manin and Katz gives the pseudo-measure of the
above conjecture when the prime $p$ is split in $K$.

The main theorem of this article (Theorem \ref{thm: main}) is the proof of the above conjecture for integers $n$ such that
the corresponding $p$-adic regulator map $r_p$ is injective.  
%
%
%
\section{Construction of the Eisenstein class}\label{section: Eisenstein class}
%
%
%

The main difficulty in the proof of the Beilinson and $p$-adic Beilinson conjectures is to construct the 
element $c^a_n \in H^1_\mot(K, M^a(n))$ for $M^a: = M(\psi^a)$ and to calculate the images $r_\infty(c^a_n)$ and $r_p(c^a_n)$
with respect to the Beilinson-Deligne and $p$-adic regulator maps.  We will use the Eisenstein symbol
as constructed by Beilinson.

We fix an integer $N \geq 3$, and let $M(N)$ be the modular curve defined over $\bbZ[1/N]$ parameterizing 
for any scheme $S$ over $\bbZ[1/N]$ the pair $(E, \levelstructure)$, where $E$ is an elliptic curve over $S$ and
$$
	\levelstructure: (\bbZ/N\bbZ)^2 \xrightarrow\cong E[N]
$$
is a full level $N$-structure on $E$, where $E[N]$ is the group of $N$-torsion points of $E$.
We let $\pr : \wt E \rightarrow M$ be the universal elliptic curve over $M$ with universal level $N$-structure
$\wt\levelstructure : (\bbZ/N\bbZ)^2\cong \wt E[N]$, and consider the motivic sheaf $\bbQ(1)$ on $\wt E$.
We let
\begin{equation}\label{eq: H}
	\sH := R^1 \pr_* \bbQ(1),
\end{equation}
and we denote by $\Sym^k \sH$ the $k$-th symmetric product of $\sH$.
Let $\varphi = \sum_{\rho \in \wt E[N] \setminus \{0\}} a_\rho [\rho]$ be a $\bbQ$-linear sum of non-zero elements in $\wt E[N]$.
For any integer $k>0$, the Eisenstein class $\Eis_\mot^{k+2}(\varphi)$ is an element
\begin{equation}
	\Eis_\mot^{k+2}(\varphi) \in H^1_\mot(M, \Sym^k \sH(1)).
\end{equation}
Although the formalism of mixed motivic sheaves or motivic cohomology with coefficients have not yet been fully developed,
one can give meaning to the above sheaves and cohomology (see \cite{BL}, \cite{BK1} for details).

Then the class $c^a_n$ may be constructed from the Eisesntein class as follows.   Let $K$ be an imaginary quadratic
field of class number one, and let $E$ be an elliptic curve defined over $\bbQ$ with complex multiplication
by the ring of integers $\cO_K$ of $K$.  We denote again by $\psi$ the Hecke character of $K$ corresponding to $E_K$
with conductor $\frf$.  We take $N \geq 3$  such that $N$ is divisible by $\frf$.
For the extension  $F: = K(E[N])$ of $K$  generated by the coordinates of the points in $E[N]$,  
we let $G_{F/K} := \Gal(F/K)$ the Galois group of $F$ over $K$.
 We fix a level $N$-structure $\levelstructure: (\bbZ/N\bbZ)^2\cong E[N]$ of $E$ over $F$, and we denote by $\levelstructure^\sigma$
 the composition of $\levelstructure$ with the action of $\sigma \in G_{F/K}$.
Then for any $\sigma \in G_{F/K}$, we denote by $\iota^{\sigma*}$ the pull-back with respect to the $F$-valued point
$\iota^\sigma: \Spec\, F \rightarrow M$ of $M$ corresponding to $(E,\levelstructure^\sigma)$.
Then the image of the sum 
$
	\iota^* :=\sum_{\sigma \in G_{F/K}} \iota^{\sigma*}
$
is invariant by the action of the 
Galois group, hence gives a pull-back morphism
$$
	H^1_\mot(M, \Sym^k \sH(1)) \xrightarrow{\iota^{*}} H^1_\mot(K, \Sym^k \iota^*\!\sH(1)).
$$
Note that on $\Spec\, K$, the motivic sheaf $\iota^*\!\sH$ is given by the motive $H^1(E)(1)$, which by definition 
corresponds to the motive $M(\psi)(1)$.  The structure of $K$-coefficients on $\iota^*\!\sH$ gives the following decomposition.


\begin{lemma}\label{lem: decomposition}
	For integers $j$ satisfying $0 \leq j \leq k/2$,  we have the decomposition of motives
	$$
		\Sym^k \iota^* \!\sH = 
			\bigoplus_{0 \leq j \leq \frac{k}{2}} M(\psi^{k-2j})\left(k-j\right),
	$$
	where we take the convention that for $k=2j$, we let $M(\psi^0)(k/2)$ be the Tate motive $\bbQ(k/2)$ with coefficients in $\bbQ$.
\end{lemma}

Let $a>0$ be an integer and we let $\frf_a$ be the conductor of $\psi^a$.  
We let $F_a := K(E[\frf_a])$ be the extension of $K$ generated by the coordinates of the points in $E[\frf_a]$, and
we let $w_{F/F_a}$ be the order of the Galois group $\Gal(F/F_a)$.
The Eisenstein classes $\Eis^{k+2}_\mot(\rho)$  are defined for points $\rho \in \wt E[N]\setminus \{0\}$ but is \textit{not}
defined for $\rho =0$.  Hence in defining $c^a_n$, we differentiate between the case when $\frf_a \neq (1)$ and $\frf_a =(1)$. 

\begin{definition}\label{def: varphi}
	We define  $\varphi_a$ as follows.
	\begin{enumerate}
		\item If $\frf_a \neq (1)$, then we fix  a primitive $\frf_a$-torsion point $\rho_a$ of $E$ and let
	$$\varphi_a := \frac{1}{w_a w_{F/F_a}} [\rho_a],$$ where we denote again by $\rho_a$ the 
	$N$-torsion point of $\wt E$ corresponding to $\rho_a$ through $\nu$ and $\wt\nu$, 
	and $w_a$ is the number of  units in $\cO_K$ which are congruent to \textit{one} modulo $\frf_a$.
		\item If $\frf_a=(1)$, then we let
	$$\varphi_a := \frac{1}{w_a w_{F/F_a}} \sum_{\rho\in\wt E[N]\setminus\{0\}}  [\rho].$$
	\end{enumerate}
\end{definition}

We define the class $c^a_n$ as follows.

\begin{definition}\label{def: can}
	For any integer $a, n$ such that $n > a >0$,  we let $k = 2n - a -2$.  Then the motive $M^a(n) := M(\psi^a)(n)$
	is a direct summand of $\Sym^k \iota^*\!\sH(1)$.
	We define the motivic class $c^a_n$ to be the image of $\Eis^{k+2}_\mot(\varphi_a)$ with respect
	to the projection
	$$
		 H^1_\mot(K, \Sym^k \iota^*\!\sH(1)) \rightarrow H^1_\mot(K, M^a(n)),
	$$
	where $\varphi_a$ is as in Definition \ref{def: varphi}.
\end{definition}

Let $p$ be a rational prime which does not divide $\frf$, and we take $N \geq 3$ to be an integer divisible by $\frf$ and prime to $p$.
In order to prove the $p$-adic Beilinson conjecture, it is necessary to calculate the images
of $c^a_n$ with respect to the Beilinson-Deligne and $p$-adic regulator maps.  The image $r_\infty(c^a_n)$
by the Beilinson-Deligne regulator map  was calculated by Deninger \cite{Den}.  We will calculate the image 
$r_p(c^a_n)$ by the $p$-adic regulator map using rigid syntomic cohomology.  

Denote by  $M^a_\cris(n)$ the crystalline realization of $M^a(n)$, which is a filtered module with a $\sigma$-linear action of
Frobenius, and let $H^1_\syn(K_\frp, M^a_\cris(n))$ be the syntomic cohomology of $K_\frp$ with coefficients in $M^a_\cris(n)$.
Then noting that $t^a_n \otimes_\bbQ \bbQ_p= M^a_\cris(n)$ in this case, there exists a canonical isomorphism
\begin{equation}\label{eq: equation one}
	t^a_n \otimes_\bbQ \bbQ_p\xrightarrow\cong H^1_\syn(K_\frp, M^a_\cris(n)).
\end{equation}
If we let  $V^a_p(n)$ be the $p$-adic \'etale realization of $M^a(n)$, then we have a canonical isomorphism 
\begin{equation}\label{eq: p-adic comparison}
	H^1_\syn(K_\frp, M^a_\cris(n)) \xrightarrow\cong H^1_f(K_\frp,V^a_p(n)),
\end{equation}
which combined with \eqref{eq: equation one} gives the exponential map \eqref{eq: exponential map}.
The syntomic regulator map
$$
	r_\syn : H^1_\mot(K, M^a(n)) \rightarrow H^1_\syn(K_\frp, M^a_\cris(n))
$$
defined by Besser (\cite{Bes1} \S7) is compatible with the $p$-adic regulator $r_p$ through the isomorphism \eqref{eq: p-adic comparison}
(\cite{Bes1}  Proposition 9.9).  Therefore, in order to calculate $\log_p \circ \, r_p(c^a_n)$, it is sufficient to calculate 
the image of $r_\syn(c^a_n)$ with respect
to \eqref{eq: equation one}.  We will calculate this image using the explicit determination of the syntomic Eisenstein class
given in \cite{BK1}.

%
%
%
\section{Eisenstein class and $p$-adic Eisenstein series}
%
%
%

In this section, we review the explicit description of the syntomic Eisenstein class in terms of $p$-adic Eisenstein
series given in \cite{BK1}.  Let $M := M(N)$ be the modular curve over $\bbZ[1/N]$ given in the previous section.
We will first describe a certain real analytic Eisenstein series $E^\infty_{k+2,l,\varphi}$. 

Let $\Gamma \subset \bbC$ be a lattice, and we denote by $A$ the area of the fundamental domain of $\Gamma$ 
divided by $\thepi := 3.14159\cdots$.  For any integer $a$ and complex number $s$ satisfying $\operatorname{Re}(s) > a/2+1$,
the Eisenstein-Kronecker-Lerch series $K^*_a(z,w,s ; \Gamma)$ to be the series
$$
	K^*_a(z,w,s ; \Gamma) :={ \sum_{\gamma \in \Gamma}}^* \frac{(\ol z+\ol\gamma)^a}{|z+\gamma|^{2s}}
	\pair{\gamma,w}
$$
where $\sum^*$ denotes the sum over $\gamma \in \Gamma$ satisfying $\gamma \not=-z$ and
 $\pair{z, w} := \exp((\ol w z - w \ol z)/A)$.  
By  \cite{We} VIII \S 12 (see \cite{BKT} Proposition 2.4 for the case $a < 0$), this series for $s$ continues meromorphically to
a function on the whole $s$-plane, holomorphic except for a simple pole at $s=1$ when $a=0$ and $w \in \Gamma$.
This function satisfies the functional equation
\begin{equation}\label{eq: functional equation}
	\Gamma(s) K^*_a(z,w,s; \Gamma) = A^{a+1-2s} \Gamma(a+1-s) K^*_a(w,z,a+1-s) \pair{w,z}.
\end{equation}

We fix a level $N$-structure $\levelstructure: (\bbZ/N\bbZ)^2 \cong \frac{1}{N}\Gamma/\Gamma$, and let
$\rho \in \frac{1}{N}\Gamma/\Gamma$.  For integers $k$ and $l$, we define the real analytic Eisenstein series 
$E^\infty_{k+2,l,\rho}$ to be the 
modular form on $M_\bbC := M(N) \otimes_\bbQ \bbC$ whose value at the test object
$(\bbC/\Gamma, dz, \levelstructure)$ is given by
\begin{equation}\label{eq: Eisenstein series}
	E^\infty_{k+2,l,\rho}(\bbC/\Gamma, dz, \levelstructure):=A^{-l}\left. \Gamma(s) K^*_{k+l+2}(0, \rho, s; \Gamma)\right|_{s=k+2}. 
\end{equation}
We let
$
	 E^\infty_{k+2,l,\varphi} := \sum_\rho a_\rho E^\infty_{k+2,l,\rho}
$
for the $\bbQ$-linear sum $\varphi = \sum_\rho a_\rho [\rho]$.  


When $l=0$, then $E^\infty_{k+2,0,\varphi}$ is a holomorphic Eisenstein series of weight 
$k+2$ on $M_\bbC$.   From the $q$-expansion, we see in this case that this Eisenstein series is
defined over $\bbQ$, and hence defines a section $E_{k+2,0,\varphi}$ in 
$\Gamma(M_\bbQ, \omega^{\otimes k} \otimes \Omega^1_{M_\bbQ})$
for $\omega := \pr_* \Omega^1_{\wt E/M}$. 
Denote by $\sH_\dR$ the de Rham realization of $\sH$, which is the coherent $\cO_{M_\bbQ}$-module
$R^1 \pr_* \Omega^\bullet_{\wt E}$  with Gauss-Manin connection 
$$
 	\nabla: \sH_\dR \rightarrow \sH_\dR \otimes \Omega^1_{M_\bbQ},
$$
and let $\Sym^k \sH_\dR$ be the $k$-th symmetric product of $\sH_\dR$ with the induced connection.
From the natural inclusion $\omega^{\otimes k} \hookrightarrow \Sym^k \sH_\dR$, we see that  $E_{k+2,0,\varphi}$
defines a section in $\Gamma(M_\bbQ, \Sym^k \sH_\dR \otimes \Omega^1_{M_\bbQ})$.

Let $p$ be a prime number not dividing $N$.   We denote by $\sH_\rig$ the filtered overconvergent 
$F$-isocrystal associated to $\sH$ on $M_{\bbZ_p}$, which is given by $\sH_\dR$ with an additional structure 
of Hodge filtration and Frobenius.  Let $H^1_\syn(M_{\bbZ_p}, \Sym^k \sH_\rig(1))$ be the rigid syntomic cohomology
of $M_{\bbZ_p}$ with coefficients in $\Sym^k \sH_\rig(1)$.  The rigid syntomic regulator is a map
$$
	r_\syn: H^1_\mot(M, \Sym^k \sH(1)) \rightarrow H^1_\syn(M_{\bbZ_p}, \Sym^k \sH_\rig(1)),
$$
and we define the syntomic Eisenstein class $\Eis^{k+2}_\syn(\varphi)$ to be the image by the syntomic regulator 
of the motivic Eisenstein class. 
We let $M^\ord_{\bbZ_p}$ be the ordinary locus in $M_{\bbZ_p}$, and 
$M^\ord_{\bbQ_p}: = M^\ord_{\bbZ_p} \otimes_{\bbZ_p} \bbQ_p$.
By \cite{BK1} Proposition A.16, a class in 
$
	H^1_\syn(M^\ord_{\bbZ_p}, \Sym^k \sH_\rig(1))
$
is given by a pair $(\cohomologyclass, \xi)$ of sections
\begin{equation}\label{eq: classes}
	\begin{split}
	\cohomologyclass &\in \Gamma(M^\ord_{\bbQ_p}, j^\dagger \Sym^k \sH_\rig(1)) \\
	\xi & \in \Gamma(M^\ord_{\bbQ_p},  \Sym^k \sH_\dR \otimes_{\bbQ_p} \Omega^1_{M^\ord}) 
	\end{split}
\end{equation}
satisfying $\nabla(\cohomologyclass) = (1 - \phi^*) \xi$.  
The $\cohomologyclass$  for the class $(\cohomologyclass, \xi)$ corresponding to the restriction to the 
ordinary locus of the syntomic Eisenstein class $\Eis^{k+2}_\syn(\varphi)$ is given as follows.


We let $p \geq 5$ be a prime not dividing $N$, and we let $\cM$ be the $p$-adic modular curve defined over 
$\bbZ_p$ parameterizing the triples
$(E_B, \eta, \levelstructure)$ consisting of an elliptic curve $E_B$ over a $p$-adic ring $B$, an isomorphism
\begin{equation}\label{eq: moduli}
	\eta:  \wh\bbG_m \cong \wh E_B
\end{equation}
of formal groups over $B$, and a level $N$-structure $\levelstructure$.
The ring of $p$-adic modular forms $V_p(\bbQ_p, \Gamma(N))$ is defined as the global section
$$
	V_p(\bbQ_p, \Gamma(N)) := \Gamma(\cM, \cO_\cM) \otimes_{\bbZ_p} \bbQ_p.
$$
The $q$-expansion gives an injection
$$
	V_p(\bbQ_p, \Gamma(N)) \hookrightarrow \bbQ_p(\zeta_N)[[q]].
$$
There exists a Frobenius action $\phi^*$ on $V_p(\bbQ_p, \Gamma(N))$ given on the $q$-expansion as 
$\phi^* = \Frob \otimes \sigma$,
where $\Frob(q) = q^p$ and $\sigma$ is the the absolute Frobenius acting on $\bbQ_p(\zeta_N)$.
The Eisenstein series $E_{k+2,0, \varphi}$ naturally defines an element in $V_p(\bbQ_p, \Gamma(N))$,
and using the fact that the differential $\partial_{\log q} := q \frac{d}{dq}$ preserves the space of
$p$-adic modular forms, we let for any integer $l \geq 0$
$$
	E_{k+l+2,l,\varphi} :=  \partial_{\log q}^l 	E_{k+2,0,\varphi}.
$$

We let
$
	E^{(p)}_{k+2,0,\varphi} := (1 - \phi^*) E_{k+2,0,\varphi}
$
and
$
	E^{(p)}_{k+l+2,l,\varphi} :=  \partial_{\log q}^l 	E^{(p)}_{k+2,0,\varphi}
$
for any integer $l \geq 0$.  Then the calculation of the $q$-expansion shows that we have
\begin{equation}\label{eq: Frobenius removal}
	E^{(p)}_{k+l+2,l,\varphi} = (1 - p^l \phi^*)  E_{k+l+2,l,\varphi}.
\end{equation}
Following the method of Katz \cite{Ka1}, we may construct a $p$-adic measure on $\bbZ_p \times \bbZ_p^\times$ 
with values in $V_p(\bbQ_p, \Gamma(N))$ satisfying the following interpolation property.

\begin{theorem}\label{theorem: Katz measure}
	There exists a $p$-adic measure $\dmu_\varphi$ on $\bbZ_p \times \bbZ_p^\times$ with values
	in $V_p(\bbQ_p, \Gamma(N))$ such that
	$$
		\int_{\bbZ_p \times \bbZ_p^\times} x^{k+1} y^l \dmu_\varphi(x,y) = E^{(p)}_{k+2,l,\varphi}
	$$
	for integers $k > 0$, $l \geq 0$.
\end{theorem}

Using this measure, we define $E^{(p)}_{k+2,l,\varphi}$ for $l<0$ as follows.

\begin{definition}[$p$-adic Eisenstein series]\label{def: p Eisenstein} 
	Let $k$ be an integer $ \geq -1$. We let
	$$
		E^{(p)}_{k+2,l,\varphi} := 	\int_{\bbZ_p \times \bbZ_p^\times} x^{k+1} y^l \dmu_\varphi(x,y)
		\,\,\in V_p(\bbQ_p, \Gamma(N)),
	$$
	where $l$ is any integer in $\bbZ$.
\end{definition}

The $p$-adic Eisenstein series satisfies the differential equation
$$
	\partial_{\log q} E^{(p)}_{k+2,l,\varphi}  = E^{(p)}_{k+3,l+1,\varphi},
$$
and the weight of $E^{(p)}_{k+2,l,\varphi}$ is $k+l+2$.
The syntomic Eisenstein class may be described using these $p$-adic Eisenstein series.
The moduli problem for $\cM$ implies that there exists a universal trivialization
$$
	\eta: \wh\bbG_m \cong \wh{\wt E}
$$
of the universal elliptic curve on $\cM$, which gives rise to a canonical section $\wt \omega$ of $\ul \omega := \pr_* \Omega^1_{\wt E/\cM}$\
corresponding to the invariant differential $d \log (1+T)$ on $\wh \bbG_m$.
Since $\cM$ is affine, there exists sections $x$ and $y$ of $\wt E$ such that the elliptic curve
$\wt E_{\bbQ_p} := \wt E \otimes \bbQ_p$ 
is given by the Weierstrass equation
$$
	\wt E_{\bbQ_p}: y^2 = 4x^3-g_2x-g_3,  \quad g_2, g_3 \in V_p(\bbQ_p, \Gamma(N))
$$
satisfying $\wt \omega = dx/y$.
Then the pull back of the $F$-isocrystal  $\sH_\rig$ to $\cM_{\bbQ_p}$  is given as
$$
	\sH_\rig = \cO_{\cM_{\bbQ_p}}\wt\omega^\vee \oplus \cO_{\cM_{\bbQ_p}}\wt u^\vee,
$$
with connection $\nabla(\wt u^\vee)  =  \wt \omega^\vee \otimes d \log q$, $\nabla(\wt \omega^\vee) = 0$, 
Frobenius $\phi^*(\wt\omega^\vee) = p^{-1} \wt\omega^\vee$, $\phi^*(\wt u^\vee) = \wt u^\vee$ and Hodge filtration
$\Fil^{-1} \sH_\rig = \sH_\rig$, $\Fil^0 \sH_\rig = \cO_{\cM_{\bbQ_p}}\wt u^\vee$, $\Fil^1 \sH_\rig = 0$
(See \cite{BK1} \S 4.3).  
If let $\wt\omega^{m,n} := \wt\omega^{\vee m} \wt u^{\vee n}$, then
the filtered $F$-isocrystal $\Sym^k \sH_\rig(1)$ on  $\cM_{\bbQ_p}$ is given by the coherent module
$$
	\Sym^k \sH_\rig(1) = \bigoplus_{j=0}^k \cO_{\cM_{\bbQ_p}}\wt\omega^{k-j, j}(1)
$$
with connection $\nabla(\wt\omega^{k-j, j}(1)) = j\,\wt\omega^{k-j +1, j-1}(1) \otimes d\log q$, Frobenius 
$$\phi^*(\wt\omega^{,k-j, j}(1)) = p^{j-k-1} \wt\omega^{k-j, j}(1),$$ 
and Hodge filtration
\begin{align*}
	\Fil^m (\Sym^k \sH_\rig(1)) = \bigoplus_{j=m+k+1}^k \cO_{\cM_{\bbQ_p}}\wt\omega^{k-j, j}(1).
\end{align*}
If we let $\wt\cohomologyclass^{k+2}_\Eis$ be the section
$$
	\wt\cohomologyclass^{k+2}_\Eis(\varphi) := \sum_{j=0}^k \frac{(-1)^{k-j}}{j!} E^{(p)}_{j+1,j-k-1,\varphi}\, \wt\omega^{k-j, j}(1),
$$
then we have 
$$
	\nabla(\wt\cohomologyclass^{k+2}_\Eis(\varphi)) =\frac{(1 - \phi^*) E_{k+2,0,\varphi}}{k!}\, \wt\omega^{0,k}(1) \otimes d\log q.
$$
The main result of \cite{BK1} is the following.

\begin{theorem}[\cite{BK1} Theorem 5.11]\label{thm: BK}
	For any integer $k>0$, the syntomic Eisenstein class 
	$$
		\Eis^{k+2}_\syn(\varphi) \in H^1_\syn(M_{\bbZ_p}, \Sym^k \sH_\rig(1))
	$$ 
	restricted to the ordinary locus $H^1_\syn(M^\ord_{\bbZ_p},  \Sym^k \sH_\rig(1))$ is represented by the pair $(\cohomologyclass, \xi)$ as in \eqref{eq: classes},
	where $\xi = E_{k+2,0,\varphi} \wt\omega^{0,k}(1)/k! \otimes d \log q$ and $\cohomologyclass$ is a section which maps to $\wt\cohomologyclass^{k+2}_{\Eis}(\varphi)$
	in $\Gamma(\cM_{\bbQ_p}, \Sym^k \sH_\rig(1))$.
\end{theorem}

The main ingredient in the proof of the above theorem is the characterization of $\xi$ 
by the residue, which by \cite{BL} 2.2.3 (see also \cite{HK1} C.1.1) and the compatibility of the Beilinson-Deligne 
regulator map with the residue morphism shows that $\xi$ represents the de Rham Eisenstein class in de Rham cohomology.
See \cite{BK1} Proposition 3.6 and Proposition 4.1 for details concerning this point.
 
%
%
%
\section{Special values of Hecke $L$-functions}\label{section: Hecke}
%
%
%

In this section, we give in Propositions 
\ref{pro: L-value} and \ref{pro: L-value-2} the precise relation between the special values of the Hecke $L$-function 
$L(\psi^a,s)$ and Eisenstein-Kronecker-Lerch series.
Assume that $K$ is an imaginary quadratic field of class number one, and let $E$ be
an elliptic curve over $\bbQ$ with good ordinary reduction at a prime $p$ with complex multiplication by 
the ring of integers $\cO_K$ of $K$.  We let $\psi$ be the Grossencharacter of $K$ associated to $E_K := E \otimes_\bbQ K$,
and we denote by $\frf$ the conductor of $\psi$.  We fix an invariant differential $\omega$ of $E$ defined over $K$.
We fix once and for all a complex embedding $\tau: K \hookrightarrow \bbC$ of the base field $K$ into $\bbC$, and we 
let $\Gamma$ be the period lattice of $E: = E \otimes_{K, \tau} \bbC$ with respect to $\omega$.  
Then we have a complex uniformization
\begin{equation}\label{eq: complex uniformization}
	 \bbC/\Gamma \xrightarrow\cong E(\bbC)
\end{equation}
such that the pull-back of the invariant differential $\omega$ coincides with $dz$. 
Note that since $E$ has complex multiplication, we have $\Gamma = \Omega \cO_K$
for some complex period $\Omega \in \bbC^\times$.

By abuse of notation, we will denote by $\psi$ and $\ol\psi$ the complex Hecke characters $\psi_\tau$ and $\ol{\psi_\tau}$
associated to $\psi$,  where $\tau$ is the fixed embedding given above.
Let $-d_K$ denote the discriminant of $K$, so that $K = \bbQ(\sqrt{-d_K})$.
The Hecke character $\psi$ is of the form $\psi(( u)) = \varepsilon( u)  u$ for any $ u \in \cO_K$
prime to $\frf$, where $\varepsilon: (\cO_K/\frf)^\times \rightarrow K^\times$ is a primitive character on $(\cO_K/\frf)^\times$.  

Let $\chi: (\cO_K/\frf_\chi)^\times \rightarrow K^\times$ be a primitive character of conductor $\frf_\chi$, and
let $f_\chi$ be a generator of $\frf_\chi$.  Then for any $u$ in $\cO_K$, 
we define the \textit{Gauss sum} $G(\chi, u)$ by
$$
	G(\chi, u) : =  \sum_{  v \in \cO_K/\frf_\chi}   \ol\chi(  v)  
	\exp\left( 2 \thepi i \Tr_{K/\bbQ}\left( u  v/f_\chi\sqrt{-d_K}\right)\right)
$$
(see \cite{L} Chapter 22 \S 1),  where we extend $\chi$ to a function on $\cO_K/\frf_\chi$ by taking $\chi(\thealpha) :=0$ for any 
$\thealpha \in \cO_K$ not prime to $\frf_\chi$.  We let $G(\chi) := G(\chi,1)$.    Then the standard fact concerning Gauss sums 
are as follows  (see for example \cite{L} Chapter 22 \S 1.) 

\begin{lemma}\label{lem: Gauss sum}
	Let the notations be as above.  
	\begin{enumerate}
		\item We have $|G(\chi)|^2= N(\frf_\chi)$.
		\item For any $u \in \cO_K$, we have
			$
				G(\chi, u) = \chi(u) G(\chi).
			$
	\end{enumerate}
\end{lemma}

As in \S \ref{section: Eisenstein class}, we let $a>0$ be an integer and $\frf_a$ be the conductor of $\psi^a$. 
Then the finite part  $\varepsilon^a$ of $\psi^a$ is a primitive character $\varepsilon^a: (\cO_K/\frf_a)^\times \rightarrow \bbC^\times$
of conductor $\frf_a$.  We fix a generator $f_a$ of $\frf_a$ and we denote by
by $G(\varepsilon^a,u)$  the corresponding Gauss sum for any $u$ in $\cO_K$.  

We let the notations be as in  \S \ref{section: Eisenstein class}.  In particular, we let $w_a$ be the number of units in $\cO_K$ 
which are congruent to \textit{one} mod $\frf_a$, and we let $w_{F/F_a}$ be the order of the Galois group $\Gal(F/F_a)$.
We again let $N \geq 3$ be a rational integer  divisible by $\frf$ and prime to $p$, 
and  $F := K(E[N])$.  We fix an isomorphism $\levelstructure: (\bbZ/N\bbZ)^2\cong\frac{1}{N}\Gamma/\Gamma$.

We first consider the case when $\frf_a \neq (1)$.  We let $\rho_a: = \Omega/f_a$ be a primitive $\frf_a$-torsion point, which 
corresponds through the uniformization \eqref{eq: complex uniformization} to a point $\rho_a \neq 0 \in E(\ol K)$.
We then have the following.

\begin{proposition}\label{pro: L-value}
	Suppose $\frf_a \neq (1)$.  Then we have
	$$
			\frac{(-1)^a}{w_a w_{F/F_a}} \sum_{\sigma \in \Gal(F/K)} E^\infty_{n,a-n, \rho_a^\sigma}(\bbC/\Gamma, dz, \nu) =	
			 \frac{G(\varepsilon^a)\ol\Omega^a}{A^{a-n}|\Omega|^{2n}} \Gamma(s) L(\psi^a,s)|_{s=n}.
	$$
\end{proposition}

\begin{proof}
	Let $w_0$ be the number of units in $\cO_K$.  By definition, we have 
	\begin{equation*}
		L(\psi^a, s) = \frac{1}{w_0} \sum_{\thealpha \in \cO_K} \frac{\psi^a(\thealpha)}{N(\thealpha)^s}
		 =  \frac{1}{w_0} \sum_{\thealpha \in \cO_K}  \frac{\varepsilon^a(\thealpha)\thealpha^a}{N(\thealpha)^s}
	\end{equation*}
	Then Lemma \ref{lem: Gauss sum} (2) gives the equality $\varepsilon^a(\thealpha) = G(\varepsilon^a, \thealpha)/G(\varepsilon^a)$.
	If we expand the definition of the Gauss sum, we see that 
	$$
		L(\psi^a, s) = \frac{1}{w_0} \sum_{\substack{\thealpha \in \cO_K\\ \thebeta\in\cO_K/\frf_a}} 
		  \frac{ \ol\varepsilon^a(\thebeta)\thealpha^a}{N(\thealpha)^s} 
		\exp\left( 2\pi i \Tr_{K/\bbQ}\left( \frac{\thealpha \thebeta}{f_a \sqrt{-d_K}}  \right) \right).
	$$
	Noting that $\cO_K$ is preserved by complex conjugation,  we see that the above is equal to
	$$
			\frac{(-1)^a}{w_0} \sum_{\thebeta \in \cO_K/\frf_a } \sum_{\thealpha\in \cO_K }
		\frac{ \ol\varepsilon^a(\thebeta)\ol \thealpha^a}{|\thealpha|^{2s}} \exp\left( \frac{2\pi}{\sqrt{d_K}} 
		\left(\frac{\thealpha\ol{\thebeta}}{\ol f_a}- \frac{\ol\thealpha \thebeta}{f_a}  		\right) \right).
	$$
	For any $\sigma \in \Gal(F/K)$, we have $\rho_a^\sigma = \rho_a^{\sigma'}$, where $\sigma'$ is the
	class of $\sigma$ in $\Gal(F_a/K)$.  If $\sigma'_\thebeta := (\thebeta, F_a/F)$ is the element
	in $\Gal(F_a/K)$ corresponding to $\thebeta \in (\cO_K/\frf_a)^\times$ through the \textit{inverse} of the Artin
	map, then by the theory of complex multiplication, we have
	$
		\rho_a^{\sigma'_{\thebeta}} = \psi(\thebeta) \rho_a.
	$
	Hence
	\begin{multline*}
		\sum_{\sigma \in \Gal(F/K)} K_a^*(0, \rho_a^\sigma, s; \Gamma)
		= w_{F/F_a} \sum_{\sigma' \in \Gal(F_a/K)} K_a^*(0, \rho_a^{\sigma'}, s; \Gamma)\\
		= w_{F/F_a} \frac{w_a}{w_0}  \sum_{\thebeta \in (\cO_K/\frf_a)^\times} \sum_{\gamma \in \Gamma} \frac{\ol\gamma^a}{|\gamma|^{2s}}
		\pair{\gamma, \psi(\thebeta)\rho_a}\\
		= w_{F/F_a} \frac{w_a}{w_0}
		  \sum_{\thebeta \in (\cO_K/\frf_a)^\times} \sum_{\gamma \in \Gamma} \frac{\ol\varepsilon^a(\thebeta)\ol\gamma^a}{|\gamma|^{2s}}
		\pair{\gamma, \thebeta \rho_a}.
	\end{multline*}
	Our assertion follows from the fact that $\Gamma = \Omega \cO_K$, $A = |\Omega|^2\sqrt{d_K}/2 \thepi$
	and the definition \eqref{eq: Eisenstein series} of the Eisenstein-Kronecker-Lerch series.
\end{proof}

The right hand side of Proposition \ref{pro: L-value} may be used to express the Hecke $L$-function on the 
other side of the functional equation as follows.

\begin{lemma}\label{lem: pre functional}
	We have
	\begin{multline}\label{eq: pre functional}
		\frac{1}{w_a w_{F/F_a}}\sum_{\sigma \in \Gal(F/K)} E^\infty_{n,a-n, \rho_a^\sigma}(\bbC/\Gamma, dz, \nu)\\=
		\frac{A^{1-n} N(\frf_a)^{a+1-n}\ol\Omega^a}{\ol f_a^a |\Omega|^{2(a+1-n)}} \left.\Gamma(s)  L(\ol\psi^a, s) \right|_{s=a+1-n}.
	\end{multline}
\end{lemma}

\begin{proof}
	We have by definition
	\begin{multline*}
		\sum_{\thebeta\in\cO_K/\frf_a} K^*_{a}(\psi(\thebeta)\rho_a, 0 , a+1-s; \Gamma) = 	\sum_{\thebeta\in\cO_K/\frf_a}
		{\sum_{\gamma \in \Gamma}}^* \frac{(\ol{\psi(\thebeta) \rho_a} + \ol\gamma)^a}{|\psi(\thebeta)\rho_a+\gamma|^{2(a+1-s)}}\\
		= \frac{ N(\frf_a)^{a+1-s} \ol\Omega^a}{\ol f_a^a |\Omega|^{2(a+1-s)}} L(\ol\psi^a, a+1-s)
	\end{multline*}
	for  $\operatorname{Re}(s) < a/2$, hence for any $s \in \mathbb{C}$ by analytic continuation.  
	Our assertion follows from the functional equation 
	$$
		\Gamma(s) K^*_{a}(0, \psi(\thebeta) \rho_a, s; \Gamma) = 
	  	A^{a+1-2s}\Gamma(a+1-s) K^*_{a}(\psi(\thebeta)\rho_a, 0 , a+1-s; \Gamma)
	$$
	and the definition \eqref{eq: Eisenstein series} of the Eisenstein-Kronecker-Lerch series.
\end{proof}

The case when $\frf_a=(1)$ is given as follows. 

\begin{proposition}\label{pro: L-value-2}
	Suppose $\frf_a = (1)$.  Then we have
	\begin{multline*}
		 \frac{1}{w_a} \sum_{\rho \in E[N] \setminus \{0\}}
		E^\infty_{n,a-n,\rho}(\bbC/\Gamma, dz, \nu)
		\\= \left( \frac{N^{a+2}}{N^{2n}} -1  \right)	 \frac{\ol\Omega^a}{A^{a-n} |\Omega|^{2n}}\Gamma(s) L(\psi^a,s)|_{s=n}.
	\end{multline*}
\end{proposition}

\begin{proof}	
	By definition,  we have
	$$
		\sum_{\rho \in \frac{1}{N} \Gamma/\Gamma} K^*_a(0, \rho, s; \Gamma) = 
		\sum_{\gamma\in \Gamma}\sum_{\rho \in \frac{1}{N} \Gamma/\Gamma} 
		\frac{\ol\gamma^a}{|\gamma|^{2s}} \pair{\gamma,\rho} \\
		=  \frac{N^{a+2}}{N^{2s}} \sum_{\gamma\in \Gamma} \frac{\gamma^a}{|\gamma|^{2s}},
	$$
	where the last equality follows from the equality
	$$
		\sum_{\rho\in\frac{1}{N}\Gamma/\Gamma} \pair{\gamma,\rho}  =
		\begin{cases}
			N^2  &  \gamma \in N\Gamma\\
			0 & \text{otherwise}
		\end{cases}
	$$
	and the fact that complex conjugation acts bijectively on  $\Gamma$.
	Our assertion follows from
	the definition \eqref{eq: Eisenstein series} of the Eisenstein-Kronecker-Lerch series.
\end{proof}

Similarly to Lemma \ref{lem: pre functional}, we have the following.

\begin{lemma}
	We have
	\begin{multline}\label{eq: pre functional-2}
		\frac{1}{w_a} \sum_{\rho \in E[N] \setminus \{0\}} E^\infty_{n,a-n,\rho}(\bbC/\Gamma, dz, \nu))\\=
		 \left( \frac{N^{a+2}}{N^{2n}} -1\right)
		\frac{A^{1-n} \ol\Omega^a}{ |\Omega|^{2(a+1-n)}} \left.\Gamma(s)  L(\ol\psi^a, s) \right|_{s=a+1-n}.
	\end{multline}
\end{lemma}

%
%
%
\section{The Main Result}
%
%
%

In this section, we give an outline of the proof of our main theorem.   We will mainly deal with the case
when $\frf_a \neq (1)$, as the case for $\frf_a = (1)$ is essentialy the same except for the factor $(N^{a+2}/N^{2n}-1)$.
 We first calculate the $p$-adic and
complex periods $\Omega_p(n)$ and $\Omega(n)$.
From the definition of $c^a_n$ and from the compatibility of the syntomic regulator with respect to pull-back morphisms, 
the restriction of the syntomic Eisenstein class through the decomposition of Lemma \ref{lem: decomposition} gives the image by the
syntomic regulator of the element $c^a_n$ in $H^1_\mot(K, M^a(n))$.  

Let the notations be as in the previous section. We denote by $\omega^*$ the class  in 
$H^1_\dR(E/\bbC)$ corresponding to $d \ol z/A$, which is in fact a class in $H^1_\dR(E/K)$.
Let $k = 2n - a -2$.   Then $\omega^{k-j+1, j+1} := \omega^{\vee k-j} \omega^{*\vee j}(1)$ for $0 \leq j \leq k$ 
form a basis of $\Sym^k \iota^*\sH_\dR(1)$. The relation between the basis $\wt\omega^{m,n}$ 
and $\omega^{m,n}$ is given by $\wt\omega^{m,n} =\Omega_\frp^{n-m}\omega^{m,n}$.  In what follows, let
$
	\varphi_a
$
be as in Definition \ref{def: varphi}.
By Theorem \ref{thm: BK}, the pull-back of the syntomic Eisenstein class 
$\Eis^{k+2}_{\syn}(\varphi_a)$ to $H^1_\syn(K_\frp, \Sym^k\iota^*\sH(1))$ is expressed by the element
$$
	\iota^* \widetilde\cohomologyclass^{k+2}_{\Eis}(\varphi_a) = \sum_{j=0}^k \frac{(-1)^{k-j}}{j!} \Omega_\frp^{2j-k} 
	E^{(p)}_{j+1, j-k-1,\varphi_a}(E, \omega, \levelstructure)\omega^{k-j+1, j+1}.
$$
Hence the element $r_\syn(c^a_n)$ in $H^1_\syn(K_\frp, M^a_\cris(n))$ corresponding by definition to the direct factor
$j = n-1$ and $j=n-a-1$ is represented by
\begin{multline*}
	  \frac{ (-1)^{n-a+1}}{\Gamma(n)} \Omega_\frp^{a} E^{(p)}_{n, a-n,\varphi_a}(E, \omega, \levelstructure) \omega^{n-a,n} \\+
	 \frac{ (-1)^{n-1}}{\Gamma(n-a)} \Omega_\frp^{-a} E^{(p)}_{n-a, -n,\varphi_a}(E, \omega, \levelstructure) \omega^{n,n-a}.
\end{multline*}
By definition of the exponential map, the element in $t^a_n \otimes \bbQ_p$ corresponding to  $r_\syn(c^a_n)$ 
through the isomorphism \eqref{eq: equation one} is
\begin{multline*}
	  \frac{ (-1)^{n-a+1}}{\Gamma(n)}\Omega_\frp^{a} E_{n, a-n,\varphi_a}(E, \omega, \levelstructure) \omega^{n-a,n} \\
	 +  \frac{ (-1)^{n-1}}{\Gamma(n-a)}\Omega_\frp^{-a} E_{n-a, -n,\varphi_a}(E, \omega, \levelstructure) \omega^{n,n-a},
\end{multline*}
where $E_{n, a-n,\varphi_a}(E, \omega, \levelstructure)$ is the element in $\widehat K^{\ur}_\frp$ satisfying 
\begin{equation}\label{eq: no p factor}
	\left( 1 - p^{a-n} \sigma^* \right)E_{n, a-n,\varphi_a}(E, \omega, \levelstructure) =  E^{(p)}_{n, a-n,\varphi_a}(E, \omega, \levelstructure).
\end{equation}
From the definition of $c^a_n$ and the discussion at the end of \S \ref{section: Eisenstein class}, this shows that we have
\begin{multline*}
	\log_p \circ \, r_p(c^a_n) =    \frac{ (-1)^{n-a+1}}{\Gamma(n)} \Omega_\frp^{a} E_{n, a-n,\varphi_a}(E, \omega, \levelstructure) \omega^{n-a,n} \\
	 + \frac{ (-1)^{n-1}}{\Gamma(n-a)}\ \Omega_\frp^{-a} E_{n-a, -n,\varphi_a}(E, \omega, \levelstructure) \omega^{n,n-a}.
\end{multline*}
This gives the following.

\begin{theorem}\label{theorem: p period}
	Let $n$ be an integer $>a$ and assume that the $p$-adic regulator $r_p$ is injective.  Then
	the $p$-adic period $\Omega_p(n)  \in  K_\frp \bigoplus K_{\frp^*}$ of the motive 
	$M^a(n)$ is given by
	\begin{multline*}
		\Omega_p(n) =   \frac{ (-1)^{n-a+1}}{\Gamma(n)}  \Omega_\frp^{a} E_{n, a-n,\varphi_a}(E, \omega, \levelstructure)  \\
	  \bigoplus\frac{ (-1)^{n-1}}{\Gamma(n-a)} \Omega_\frp^{-a} E_{n-a, -n,\varphi_a}(E, \omega, \levelstructure).
	\end{multline*}
\end{theorem}

\begin{proof}
	The theorem follows from the definition given in \eqref{eq: p-adic period} of the $p$-adic period, noting that 
	we have an isomorphism
	$$
		(K \otimes_\bbQ \bbQ_p)\, \omega^a_{\tg,n} \cong \bbQ_p \, \omega^{n,n-a} \bigoplus \bbQ_p\, \omega^{n-a,n}
	$$
	induced from the canonical splitting $K \otimes \bbQ_p = K_\frp \bigoplus K_{\frp^*} \cong \bbQ_p \bigoplus \bbQ_p$.
\end{proof}

The calculation of the complex period, originally due to Deninger \cite{Den} may be done in a similar fashion.
If we let $\Gamma =  \Omega\cO_K$ as in the previous section, then the Betti homology of $E$ is given by
$
	H_1^B(E(\mathbb{C}), \bbZ) = \Gamma.
$
We let $\gamma_1 := \Omega \in \Gamma$, which is a generator of $\Gamma$ as a $\cO_K$-module.  If we fix a $\thebeta \in \cO_K$
such that $\cO_K := \bbZ \oplus \bbZ\thebeta$ and if we let $\gamma_2: = \thebeta(\gamma_1)$ where $\thebeta$ acts through 
the complex multiplication of $E$, then we have $\Gamma := \bbZ \gamma_1 \oplus \bbZ \gamma_2$
as a $\bbZ$-module.  The period relation gives the equality
$$
	\begin{pmatrix}
		\gamma_1 \\ \gamma_2 
	\end{pmatrix}
	=
	\begin{pmatrix}
		\Omega & \ol\Omega/A\\ 
		\tau\Omega & \ol\tau\ol\Omega/A 
	\end{pmatrix}
	\begin{pmatrix}
		\omega^\vee \\ \omega^{*\vee}
	\end{pmatrix}.
$$
The $K$-basis $\gamma_1$ induces a $K$-basis of $M^a_B(n) \subset \Sym^k H^1_B(E(\bbC), \bbQ(1))$, which we 
denote by $\omega_B$.  Then the inclusion
$$
	M^a_B(n) \otimes_\bbQ \bbR \hookrightarrow  t^a_n \otimes_\bbQ \bbR 
$$
maps $\omega_B$ to
$$
	 \Omega^{n-a}  \frac{\ol\Omega^n}{A^n} \omega^{n-a,n} + \Omega^n \frac{\ol\Omega^{n-a}}{A^{n-a}} \omega^{n,n-a}. 
$$
Furthermore, one may prove the following.

\begin{theorem}[Deninger \cite{Den}]\label{thm: Deninger}
	The image $ r_\infty(c^a_n) $ in 
	$$
		H^1_\sD(K \otimes_\bbQ \bbR, V^a_\infty(n)) \cong (t^a_n \otimes_\bbQ \bbR) / (M^a_B(n) \otimes_\bbQ \bbR)
	$$ 
	of $c^a_n$ by the Beilinson regulator \eqref{eq: Beilinson regulator} is represented by the element
	\begin{multline*}
		\frac{(-1)^{n-1}}{\Gamma(n)} A^{n-a} E^\infty_{n,a-n,\varphi_a}(\bbC/\Gamma, dz, \levelstructure)\, \omega^{n-a,n}\\+
		\frac{(-1)^{n-a+1}}{\Gamma(n-a)} A^n E^\infty_{n-a,-n,\varphi_a}(\bbC/\Gamma, dz, \levelstructure)
		 \, \omega^{n,n-a}
	\end{multline*}
	in $t^a_n \otimes_\bbQ \bbC$.
\end{theorem}

\begin{proof}
	The Eisenstein class in this paper defined using the elliptic polylogarithm is related to the Eisenstein class defined by Beilinson and Deninger.
	The theorem is then a special case of the weak Beilinson conjecture for Hecke character associated to imaginary quadratic fields
	proved by Deninger \cite{Den} (see also \cite{DW} for the case of an elliptic curve  defined over $\bbQ$ with complex multiplication.)  
	The theorem may also be proved by explicitly calculating the Hodge realization of the elliptic polylogarithm
	\cite{BL} (see also \cite{BKT} Theorem A 29.)
\end{proof}

By taking the determinant of the complex \eqref{eq: complex determinant} with respect to the basis 
$r_\infty(c^a_n)$, $\omega_B$, $\omega^{n-a,n}$ and $\omega^{n,n-a}$, the above calculation
and the definition of the complex period give the following.  

\begin{proposition}\label{prop: period}
	The complex period $\Omega_\infty(n)$  of $M^a(n)$ in $K \otimes_\bbQ \bbR$ is given by
	$$
		\Omega_\infty(n) =   (-1)^{a+n-1} G(\varepsilon^a) L(\psi^a,n) \, \bigoplus \,(-1)^{n+1}  G(\ol\varepsilon^a) L(\ol\psi^a,n )
	$$
	in $K \otimes_\bbQ \bbC \cong \bbC \bigoplus \bbC$ if $\frf_a \neq (1)$.   A similar formula holds for the case when $\frf_a = (1)$, 
	but with a factor $(N^{a+2}/N^{2n}-1)$ multiplied to the $L$-value.
\end{proposition}

\begin{proof}
	The assertion follows from Theorem \ref{thm: Deninger} by explicit calcuation, using
	the definition  of $\varphi_a$ (Definition \ref{def: varphi}), the calculation of the complex period above,
	and the relation
	between  Eisenstein-Kronecker-Lerch series and special values of $L$-functions
	(Proposition \ref{pro: L-value} if $\frf_a\neq(1)$, or Proposition \ref{pro: L-value-2} if $\frf_a = (1)$). 
\end{proof}

This gives the following corollary, which we stated in Theorem \ref{thm: Deninger one}.

\begin{corollary}\label{cor: Deninger}
	If $\frf_a \neq (1) $, then we have
	$$
		\frac{L(\psi^a,n)}{\Omega_\infty(n)} =\frac{ (-1)^{a+n-1}}{G(\varepsilon^a)} \,\in\, K^\times \subset \,K \otimes_\bbQ \bbC.
	$$
	A similar formula holds for the case when $\frf_a = (1)$, 
	but with multiplication by $(N^{a+2}/N^{2n}-1)^{-1}$ on the right hand side.
\end{corollary}

\begin{proof}
	The equalitiy follow from the calculation of the complex period in Proposition \ref{prop: period}.
	Since $\varepsilon^a$ is a primitive Hecke character
	with values in $K$, we see that this value is in $K$.
\end{proof}

We next construct the $p$-adic measure $\dmu^a$ which appears in the formulation of Conjecture \ref{conj: PR}.
Since $E$ has good ordinary reduction at $p$, the prime $p$ splits as $p =\frp \frp^*$ in $K$.
In what follows, we fix once and for all complex and $p$-adic embeddings of our coefficient $K$ as follows.  We let $\tau: K \hookrightarrow \bbC$ as in 
\S \ref{section: Hecke} and an embedding $K \hookrightarrow \bbC_p$ mapping $\frp$ to a prime in $\bbC_p$.  With this convention,
we may regard the complex and $p$-adic periods as elements respectively in $\bbC$ and $\bbC_p$, by taking the first
components in Proposition \ref{prop: period} and Theorem \ref{theorem: p period}.

Let $\wh E$ be the formal group of $E$ over $\cO_K$, and let $\wh K_\frp^\ur$ be the $p$-adic completion of the maximal 
unramified extension $K^\ur_\frp$ of $K_\frp$, which we regard as a subfield of $\bbC_p$ through our fixed embedding.
Since $p$ is an ordinary prime, there exists an isomorphism of formal groups $\eta$ over 
$\cO_{\wh K_\frp^\ur}$ 
$$
	\eta:  \wh E \xrightarrow\cong  \wh\bbG_m 
$$
given by a power series $\eta(t) = \exp(\lambda(t)/\Omega_\frp)-1$, 
where $\Omega_\frp$ is a $p$-adic period of $E$ which is an element in 
$\cO_{\wh K_\frp^\ur}$ satisfying
\begin{equation}\label{eq: Frobenius action}
	\Omega_\frp^{\sigma-1} = \psi(\frp)p^{-1}.
\end{equation}
The above isomorphism gives the equality
\begin{equation}\label{eq: normalization}
	\eta^*\left(d \log(1+T)\right) =\omega/\Omega_\frp.
\end{equation}

Again let $N \geq 3$ be an integer as in \S \ref{section: Eisenstein class}  divisible by $\frf$ and prime to $p$.
By \cite{Ka1} 5.10.1, the value of the $p$-adic Eisenstein series  $E_{k+2,l,\varphi_a}$ 
at the test object $(E, \omega, \levelstructure)$ is defined by
\begin{equation}\label{eq: value of specialization}
	E_{k+2,l,\varphi_a}(E, \omega, \levelstructure) := \Omega_\frp^{k+l+2} E_{k+2,l,\varphi_a}(E, \eta, \levelstructure).
\end{equation}
In addition, the comparison theorem \cite{Ka1} 8.0.9 states that this value for integers $k>0, l\geq 0$ is an element in $F:= K(\frf)$ satisfying 
the equality
$$
	E_{k+2,l,\varphi_a}(E, \omega, \levelstructure) = E^\infty_{k+2,l,\varphi_a}(E, \omega, \levelstructure).
$$
Then the calculation above and Theorem \ref{theorem: Katz measure} gives the following.

\begin{proposition}
	We let $\varphi_a$ be as in Definition \ref{def: varphi}, and 
	we denote again by $\dmu_{\varphi_a}$ the $p$-adic measure 
	on $\bbZ_p \times \bbZ_p^\times$ obtained as the value of $\dmu_{\varphi_a}$ of 
	Theorem \ref{theorem: Katz measure} at $(E, \omega, \levelstructure)$.  
	If $\frf_a\neq (1)$, then we have
	$$
		\frac{(-1)^a}{\Omega_\frp^a} \int_{\bbZ_p \times \bbZ_p^\times} x^{n-1} y^{a-n} \dmu_{\varphi_a}(x,y) = 
		G(\varepsilon^a)\left(1 - \frac{\psi(\frp)^a}{p^n} \right) 
		\frac{\ol\Omega^a \Gamma(n)L(\psi^a, n)}{ A^{a-n} |\Omega|^{2n} } 
	$$
	for integers $a \geq n > 0$.
	A similar formula holds for the case when $\frf_a = (1)$, 
	but with multiplication by $(N^{a+2}/N^{2n}-1)$ on the right hand side.
\end{proposition}

\begin{proof}
	The relation between the action of the Frobenius on $V_p(\bbQ_p, \Gamma(N))$ and its specialization is given by
	$$
		\phi^*(E_{k+2,l,\varphi} )(E, \eta, \levelstructure) = E_{k+2,l,\varphi} (E,\eta, \levelstructure^{\sigma_\frp}),
	$$
	since $(E,\omega)$ is define over $K$ and hence $E^{\sigma_\frp} = E$ and $\omega^{\sigma_\frp} = \omega$.
	Then from the definition of the specialization of $p$-adic modular forms \eqref{eq: value of specialization} and the action of the Frobenius
	on the $p$-adic period  \eqref{eq: Frobenius action}, we have
	\begin{multline*}
		p^l \phi^*(E_{k+2,l,\varphi} )(E, \omega, \levelstructure) 
		= p^l \left(\Omega_\frp^{\sigma_\frp}\right)^{k+l+2} E_{k+2,l,\varphi} (E, \eta, \levelstructure^{\sigma_\frp})\\
		= \psi(\frp)^{k+l+2} p^{-k-2} \Omega_\frp^{k+l+2} E_{k+2,l,\varphi^{\sigma_\frp}} (E, \eta, \levelstructure).
	\end{multline*}
	Applying the above calculation to the case $a = k+l+2$ and $n= k+2$,
	our assertion now follows from Theorem \ref{theorem: Katz measure}, noting  the definition of the $p$-adic Eisenstein series 
	\eqref{eq: Frobenius removal}, the definition of $\varphi_a$, and  the fact that the sum over all $\sigma \in G_{F/K}$ of $\varphi_a^\sigma$ 
	is invariant by the action
	of $\sigma_\frp$.
\end{proof}

\begin{proposition}\label{pro: measure}
	Let $\bbZ_p^\times \times \bbZ_p^\times \rightarrow \bbZ_p^\times$ be the surjection defined by $\rho: (x,y) \mapsto x/y$.
	We define the measure $\dmud^a$ on $\bbZ^\times_p$ by
	$$
		\int_{\bbZ_p^\times} f(w) \dmud^a(w) = 
		\frac{(-1)^a}{G(\varepsilon^a)} \int_{\bbZ^\times_p \times \bbZ_p^\times} x^{-1} y^a \rho^*(f)(x, y) \dmu_{\varphi_a}(x,y).
	$$
	If $\frf_a \neq (1)$, then this measure satisfies the interpolation property 
	$$	
		\frac{1}{\Omega^a_\frp} \int_{\bbZ_p^\times} w^n \dmud^a(w) = 
		\left(1 - \frac{\psi(\frp)^a}{p^n} \right)
		 \left(1 - \frac{\ol\psi(\frp^{*})^a}{p^{a+1-n}} \right)
		   \frac{\ol\Omega^a \Gamma(n) L(\psi^a, n)}{A^{a-n}  |\Omega|^{2n} }.
	$$
	for any integers $a$ and $n$ such that $a \geq n > 0$.  
	A similar formula holds for the case when $\frf_a = (1)$, 
	but with multiplication by $(N^{a+2}/N^{2n}-1)$ on the right hand side.
\end{proposition}

\begin{proof}
	The equality is obtained from the definition of $\dmud^a$ and in
	calculating the restriction of the measure $\dmu_{\varphi_a}$ to $\bbZ^\times_p \times \bbZ^\times_p$.
	The calculation directly follows from Katz \cite{Ka1} 8.7.6, using the functional equation (see Remark \ref{rem: functional equation} below.)	
	One may also do the calculation using an alternative construction of Katz $p$-adic measure (\cite{BKo} Proposition 3.5 and Theorem 3.7),
	again after using the functional equation.
\end{proof}

\begin{remark}\label{rem: functional equation}
	Combining Proposition \ref{pro: L-value} and \eqref{eq: pre functional} (or if $\frf_a=(1)$, then 
	Proposition \ref{pro: L-value-2} and \eqref{eq: pre functional-2}), we obtain the functional equation
	$$
		\frac{\Gamma(n) L(\psi^a,n)}{A^{a-n} |\Omega|^{2n} } 
		= \frac{ N(\frf_a)^{a+1-n} \Gamma(a+1-n) L(\ol\psi^a, a+1-n) }{(-1)^aG(\varepsilon^a)\ol f_a^aA^{n-1}|\Omega|^{2(a+1-n)}}.
	$$
	We regard $f_a$ and $\ol f_a$ in $\cO_K$ as elements in $\bbZ_p^\times$ through the canonical isomorphism
	$\cO_{K_\frp} \cong \bbZ_p$.
	Denote by $\dwtmu_{\varphi_a}$ the measure on $\bbZ_p^\times \times \bbZ_p^\times$ 
	obtained as the pull-back of $f_a^{-1} \dmu_{\varphi_a}$ through the isomorphism
	$$
		\bbZ_p^\times \times \bbZ_p^\times  \xrightarrow\cong 	\bbZ_p^\times \times \bbZ_p^\times
	$$
	given by $(x,y) \mapsto (\ol f_a x,f_a^{-1} y)$.  Since $A = |\Omega|^2 \sqrt{d_K}/2 \thepi$, 
	if we let  $k_1 = a+1-n$ and $k_2 = 1-n$, then the interpolation property of $\dwtmu_{\varphi_a}$ at $(E, \omega, \levelstructure)$ becomes
	\begin{multline*}
		\frac{1}{\Omega_\frp^{k_1-k_2}} \int_{\bbZ_p^\times \times \bbZ_p^\times} x^{-k_2} y^{k_1-1} 
		\dwtmu_{\varphi_a}(x,y) \\= 
		  \left(1 - \frac{\psi(\frp)^{k_1-k_2}}{p^{1-k_2}} \right) 
		 \left(1 - \frac{\ol{\psi}(\frp^{*})^{k_1-k_2}}{p^{k_1}} \right) \left( \frac{\sqrt{d_K}}{2\thepi}\right)^{k_2}
		 \frac{\Gamma(k_1) L(\ol\psi^{{k_1-k_2}}, k_1)}{\Omega^{k_1-k_2}}
	\end{multline*}
	for $k_1 > -k_2 \geq 0$. This coincides with the interpolation property of the two-variable $p$-adic measure 
	constructed by Katz and Yager  (see \cite{Ya} \S 1.)
\end{remark}

If $\frf_a \neq (1)$, then the measure $\dmud^a$ defined in Proposition \ref{pro: measure} satisfies the condition of Conjecture \ref{conj: PR}.
If $\frf_a = (1)$, then we need to cancel the factor $(N^{a+2}/N^{2n}-1)$ which appears in the interpolation formula.

\begin{definition}\label{def: measure}
	We define the pseudo-measure $\dmu_a$ on $\bbZ_p^\times$ as follows.
	\begin{enumerate}
		\item If $\frf_a \neq (1)$, then we let $\dmu^a := \dmud^a$.
		\item If $\frf_a = (1)$, then we let $\dmu^a_N$ be the measure on $\bbZ_p^\times$ defined by
		$$
			\int_{\bbZ_p^\times} x^n \dmu^a_N = \left(\frac{N^{a+2}}{N^{2n}}-1 \right)
		$$
		for any integer $n$.  We define $\dmu^a$ to be the pseudo-measure  on $\bbZ_p^\times$ 
		obtained as the quotient of $\dmud^a$ by $\dmu^a_N$ (see for example \cite{Col} \S 1.2 for the definition of 
		a pseudo-measure).
	\end{enumerate}
	When $\frf_a \neq (1)$, then $\dmu_a$ is by definition a $p$-adic measure on $\bbZ_p^\times$.
\end{definition}

We now have the following.

\begin{theorem}\label{thm: main}
	Let $a>0$ be an integer and let $\dmu^a$ be the pseudo-measure on $\bbZ_p^\times$ defined in Definition \ref{def: measure}.  
	If we let 
	$$
		L_p(\psi^a \otimes \chi_\cyc^n) := \int_{\bbZ_p^\times} w^n \dmu^a(w),
	$$
	then we have
	$$
		\frac{L_p(\psi^a \otimes \chi_\cyc^n)}{\Omega_p(n)}=
			\left(1 - \frac{\psi(\frp)^a}{p^n} \right) \left(1 - \frac{\ol\psi(\frp^{*})^a}{p^{a+1-n}} \right)  \frac{\Gamma(n) L(\psi^a,n)}{\Omega_\infty(n)}
	$$
	for any integer $n > a$ such that the $p$-adic regulator map $r_p$ is injective.
\end{theorem}

\begin{proof}
	We first consider the case when $\frf_a\neq(1)$. By definition of the $p$-adic Eisenstein series in Definition \ref{def: p Eisenstein}, the moments of the measure $\dmu_{\varphi_a}$ constructed 
	in Theorem \ref{theorem: Katz measure} is given by $p$-adic Eisenstein series.   As in \eqref{eq: no p factor},
	let $E_{n,a-n,\varphi_a}(E, \omega, \levelstructure)$ be the element in $\wh K_\frp^\ur$ satisfying
	\begin{equation*}
		\left( 1 - p^{a-n} \sigma^* \right)E_{n, a-n,\varphi_a}(E, \omega, \levelstructure) =  E^{(p)}_{n, a-n,\varphi_a}(E, \omega, \levelstructure).
	\end{equation*}
	Then the compatibility between the Frobenius on the modular curve and a point, as well as the restriction of the
	measure on $\bbZ_p \times \bbZ_p^\times$ to $\bbZ_p^\times\times\bbZ_p^\times$ shows that we have the relation
	\begin{multline*}
		\frac{1}{\Omega^a_\frp} \int_{\bbZ_p^\times} w^n \dmu^a(w)
		\\ =  \frac{1}{G(\varepsilon^a)}
		\left(1 - \frac{\psi(\frp)^a}{p^n} \right) \left(1 - \frac{\ol\psi(\frp^{*})^a}{p^{a+1-n}} \right) E_{n,a-n,\varphi_a}(E, \omega, \levelstructure).
	\end{multline*}
	The calculation of the $p$-adic period in Theorem \ref{theorem: p period} shows that
	$$
		\Omega_p(n) =\frac{ (-1)^{a-n+1}}{\Gamma(n)} \Omega^a_\frp E_{n,a-n,\varphi_a}(E, \omega, \levelstructure)
	$$
	for our fixed embedding $K \hookrightarrow \bbC_p$.   This proves in particular that
	$$
	\frac{(-1)^{a-n+1}}{\Omega_p(n)} \int_{\bbZ_p^\times} w^n \dmu^a(w) =
		\left(1 - \frac{\psi(\frp)^a}{p^n} \right) \left(1 - \frac{\ol\psi(\frp^{*})^a}{p^{a+1-n}} \right) \frac{\Gamma(n)}{G(\varepsilon^a)}.
	$$
	Our assertion now follows from Corollary \ref{cor: Deninger}.
	The case for $\frf_a=(1)$ follows in a similar fashion, noting the interpolation property of $\dmud^a$ and $\dmu^a_N$.
\end{proof}

%
%


\end{document}